\documentclass[11pt, leqno]{article}
\usepackage{amsmath, amsfonts, amsthm}
\usepackage{xcolor}
 \numberwithin{equation}{section}

\begin{document}

\author{Ajai Choudhry, Iliya Bluskov and Alexander James}
\title{A quartic diophantine equation\\ inspired by Brahmagupta's identity}
\date{}
\maketitle

\begin{abstract}
In this paper we obtain several parametric solutions of the quartic diophantine equation
$(x_1^4+x_2^4)(y_1^4+y_2^4)=z_1^4+z_2^4$. We also show how infinitely many parametric solutions of this equation may be obtained by using elliptic curves. 
\end{abstract}


\section{Introduction}

If $x_1,\,x_2,\, y_1,\,y_2$ are any four arbitrary integers, it is easy to find integers $z_1$ and $z_2$ such that
\begin{equation}
(x_1^2+x_2^2)(y_1^2+y_2^2)=z_1^2+z_2^2. \label{quadeqn}
\end{equation}
This follows from the well-known identity,
\begin{equation}
(x_1^2+x_2^2)(y_1^2+y_2^2)=(x_1y_1+x_2y_2)^2+(x_1y_2-x_2y_1)^2, \label{Brahmident}
\end{equation}
given by Brahmagupta in the seventh century AD. 

A natural question arises whether there exist any integers satisfying the relations,
\begin{equation}
(x_1^n+x_2^n)(y_1^n+y_2^n)=z_1^n+z_2^n, \label{geneqn}
\end{equation}
when the  exponent $n > 2$. 

It is readily seen that if  $x_1,\,x_2,\, y_1,\,y_2$ are any four arbitrary integers,  there do not exist integers $z_1$ and $z_2$ such that 
\begin{equation}
(x_1^3+x_2^3)(y_1^3+y_2^3)=z_1^3+z_2^3. \label{cubeqn}
\end{equation}
For instance, if each of the integers $x_1,\,x_2,\, y_1,\,y_2$ is $ \equiv 1 (\rm mod \; 3)$, then the left-hand side of \eqref{cubeqn} is $ \equiv 4 (\rm mod \; 9)$ while for arbitrary integers $z_1$ and $z_2$, the right-hand side of \eqref{cubeqn} can never be equal to $ 4 (\rm mod \; 9)$. Thus,  we cannot have an identity \eqref{cubeqn} with $x_1,\,x_2,\, y_1,\,y_2$ being  four arbitrary integers. 

It is, however, easy to find  integers $x_i,\, y_i,\,z_i,\; i=1,\,2$ such that Eq.~\eqref{cubeqn} is satisfied. In fact, given any two  arbitrary integers $y_1,\,y_2$, we can find integers  $x_i,\,z_i,\; i=1,\,2$ such that Eq.~\eqref{cubeqn} is satisfied. The complete solution of the cubic diophantine equation,
\begin{equation}
a(x^3+y^3)=b(u^3+v^3)
\end{equation}
has been given by Choudhry \cite{Ch}. Using this solution, and taking $a=y_1^3+y_2^3,\; b=1$, we get, after suitable renaming of the variables, the complete solution of Eq.~\eqref{cubeqn} with $y_1$ and $y_2$ being arbitrary integers. Thus, the existing results immediately give a complete solution of Eq.~\eqref{cubeqn}.

This paper is concerned with Eq.~\eqref{geneqn} when $n=4$. It seems that this equation has not been considered in the existing literature. We will obtain several parametric solutions of this equation in Section 2. In fact, by using elliptic curves, we will show  that infinitely many parametric solutions of this equation can be obtained. We conclude the paper with some related open problems suggested by Brahmagupta's identity.

\section{The quartic diophantine equation $(x_1^4+x_2^4)(y_1^4+y_2^4)=z_1^4+z_2^4$}

We will now consider the quartic diophantine equation,
\begin{equation}
(x_1^4+x_2^4)(y_1^4+y_2^4)=z_1^4+z_2^4. \label{quarteqn}
\end{equation}

We first observe that, as in the case of Eq.~\eqref{cubeqn}, if  $x_1,\,x_2,\, y_1,\,y_2$ are four arbitrary integers,  there do not exist integers $z_1$ and $z_2$ such that Eq.~\eqref{quarteqn} is satisfied. For instance, if all  the integers $x_1,\,x_2,\, y_1,\,y_2$ are odd, then the left-hand side of \eqref{cubeqn} is $ \equiv 4 (\rm mod \; 16)$ while for arbitrary integers $z_1$ and $z_2$, the right-hand side of \eqref{cubeqn} can never be equal to $ 4 (\rm mod \; 16)$.

We further  note that if $(x_1,\,x_2,\,y_1,\,y_2,\,z_1,\,z_2)=(\alpha_1,\,\alpha_2,\,\beta_1,\,\beta_2,\,\gamma_1,\,\gamma_2)$ is any solution of Eq.~\eqref{quarteqn}, then taking $k_1$ and $k_2$ as any two  arbitrary nonzero integers, $(x_1,\,x_2,\,y_1,\,y_2,\,z_1,\,z_2)=(k_1\alpha_1,\,k_1\alpha_2,\,k_2\beta_1,\,k_2\beta_2,\,k_1k_2\gamma_1,$ $k_1k_2\gamma_2)$  is also a solution of Eq.~\eqref{quarteqn}. It follows that  integer solutions of Eq.~\eqref{quarteqn} may be obtained from any rational solution of \eqref{quarteqn} by suitably choosing the integers $k_1$ and $k_2$.

Since the left-hand side of Eq.~\eqref{quarteqn} may be expressed as a sum of four biquadrates, each integer solution of Eq.~\eqref{quarteqn} immediately yields a solution of the diophantine equation,
\begin{equation}
u_1^4+u_2^4+u_3^4+u_4^4=z_1^4+z_2^4. \label{quarteqnH}
\end{equation}
Haldeman has described a method of obtaining two parametric solutions of Eq.~\eqref{quarteqnH} (as mentioned by Dickson \cite[p.\ 653]{Di}) but these parametric solutions do not readily yield any nontrivial solution of Eq.~\eqref{quarteqn}.

We also note that integer solutions of Eq.~\eqref{quarteqn} with $x_1=x_2=1$ have been obtained by Izadi and Nabardi \cite{Iz}. We will accordingly consider solutions of Eq.~\eqref{quarteqn} in which $x_1 \neq x_2$ and $y_1 \neq y_2$. A quick computer search yielded several solutions in small integers. In Table I below, we give 12  integer solutions of Eq.~\eqref{quarteqn} found by computer trials.

\begin{center}
\begin{tabular}{|c|c|c|c|c|c|c|} \hline
\multicolumn{7}{|c|}{\bf  Table I: Solutions of Eq.~\eqref{quarteqn}} \\
\hline
S.~No. & $x_1$ &$x_2$ &$y_1$ &$y_2$ &$z_1$ &$z_2$\\
\hline
1& 1 & 2 & 5 & 6 & 8 & 13 \\
2& 1  & 2 & 25 & 28 & 39 & 62\\
3& 1  & 4  & 4 & 15  &   49 & 52  \\
4& 1  & 5  & 16 &  29 &  97 & 141  \\
5 & 1  & 8 & 65 &  264 & 448 & 2113  \\
6 & 1  & 10 & 8  & 11  &  2 & 17  \\
7 & 2  & 5  &  16  &  19  &  78 & 97  \\
8 & 3 & 5 & 17 & 28& 13 & 149 \\
9 & 3 & 10 & 6 & 17 & 8 & 171 \\
10 & 3 & 14 & 5 & 6 & 39 & 92 \\
11 & 5 & 6 & 6 & 13 & 16 & 87 \\
12 & 8 & 11 & 13 & 15 & 163 & 167 \\
\hline
\end{tabular}
\end{center}

The numerical solutions  of Eq.~\eqref{quarteqn} obtained by computer search suggested the existence of parametric solutions. Accordingly, in the next two subsections, we explore parametric solutions of Eq.~\eqref{quarteqn}. We   obtain several  parametric solutions, and show how infinitely many parametric solutions of Eq.~\eqref{quarteqn} may be obtained.

\subsection{A family of parametric solutions}
Using Brahmagupta's identity, we may  write,
 \begin{equation}
(x_1^4+x_2^4)(y_1^4+y_2^4)=(x_1^2y_1^2+x_2^2y_2^2)^2+(x_1^2y_2^2-x_2^2y_1^2)^2. \label{quartident}
\end{equation}
To obtain a solution of Eq.~\eqref{quarteqn}, we must choose integers  $x_1,\,x_2,\, y_1,\,y_2$  such that the following equations are simultaneously satisfied:
\begin{align}
x_1^2y_1^2+x_2^2y_2^2&=z_1^2, \label{eqz1}\\
x_1^2y_2^2-x_2^2y_1^2&=z_2^2. \label{eqz2}
\end{align}
Thus, $(x_1y_1,\,x_2y_2,\,z_1)$ is a Pythagorean triple, and we may write,
\begin{equation}
x_1y_1=p^2-q^2, \quad x_2y_2=2pq,\quad z_1=p^2+q^2,
\end{equation}
where $p$ and $q$ are arbitrary  parameters. It follows that
\begin{equation}
x_1=f(p-q), \quad x_2=2gq,\quad y_1=(p+q)/f, \quad y_2=p/g, \label{valxy}
\end{equation}
where $f$ and $g$ are some rational numbers. 

On substituting the values of $x_i, \,y_i$ given by \eqref{valxy} in Eq.~\eqref{eqz2}, we get,
\begin{multline}
f^4p^4 - 2f^4p^3q + (f^2 + 2g^2)(f^2 - 2g^2)p^2q^2 \\
- 8g^4pq^3 - 4g^4q^4 = f^2g^2z_2^2, \quad \quad \quad \quad \quad \label{quarticec1}
\end{multline}
and on making the transformation,
\begin{equation}
p=qU,\quad g=fm,\quad z_2=q^2V/m, \label{tr1}
\end{equation}
we get the equation,
\begin{equation}
V^2=U^4 - 2U^3 - (2m^2 + 1)(2m^2 - 1)U^2 - 8m^4U - 4m^4, \label{quarticec2}
\end{equation}
where $m$ is some rational number.

Now Eq.~\eqref{quarticec2} represents a quartic model of an elliptic curve over the field $\mathbb{Q}(m)$, and the birational transformation defined by
\begin{equation}
\begin{aligned}
U& =(X + Y+8m^4)/(2X-8m^4), \\
V&=\{X^3 - 12m^4X^2 + 8m^4(4m^4 - 5)X \\
& \quad \quad - 24m^4Y - 128m^8\}/\{4(X-4m^4)\}^2,
\end{aligned}
\label{biratUV}
\end{equation}
and 
\begin{equation}
\begin{aligned}
X&=2U^2 - 2U + 2V,\\
 Y& = 4U^3 - 6U^2 + 4UV - 2(2m^2 + 1)(2m^2 - 1)U - 2V - 8m^4,
\end{aligned}
\label{biratXY}
\end{equation}
reduces the quartic curve \eqref{quarticec2} to the Weierstrass model,
\begin{equation}
Y^2 =X^3 - (2m^2 + 1)(2m^2 - 1)X^2 + 32m^4X. \label{ecweier} 
\end{equation}

It is readily seen that a rational point  $P$ on the elliptic curve \eqref{ecweier} is given by
\[
X = 4(m^4 - 2)^2/9, Y = 4(m^4 - 2)(2m^8 - 17m^4 - 10)/27.
\]

When $m=1$, the curve \eqref{ecweier} reduces to
\begin{equation}
Y^2 = X^3 - 3X^2 + 32X, \label{ecspl}
\end{equation}
and a rational point on the  curve \eqref{ecspl}  corresponding to the point $P$ is $(4/9,\, 100/27)$. Since this rational point on the elliptic curve \eqref{ecspl} does not have integer coordinates,  it follows from the Nagell-Lutz theorem \cite[p. 56]{Si} on elliptic curves that the point $(4/9,\, 100/27)$ is not a point of finite order. In fact, the rank of the elliptic curve \eqref{ecspl}, as determined by the software SAGE,  is 1. We can thus find infinitely many rational points on the curve \eqref{ecspl} using the group law.

Since in the special case $m=1$, the point on the curve \eqref{ecspl} corresponding to the point $P$ is not of finite order, it follows that the point $P$ on the curve \eqref{ecweier} cannot  be a point of finite order. We can thus generate infinitely many rational points on the curve \eqref{ecweier} using the group law. Each of these rational points will yield a corresponding rational point on the curve \eqref{quarticec2} and using the relations \eqref{biratUV}, \eqref{tr1} and \eqref{valxy} (in that order), we can find infinitely many parametric solutions of the diophantine equation \eqref{quarteqn}.

The  parametric solution of Eq.~\eqref{quarteqn}, obtained from the  point $P$ on the curve \eqref{ecweier},  may be written,  after removing certain common factors and making appropriate changes of sign,  as follows:
\begin{equation}
\begin{aligned}
x_1& = 6m,  \quad  &x_2& =(m^2 - 2m + 2)(m^2 + 2m + 2), \\
y_1 & = m(m^4 - 2), \quad &y_2& = m^4 + 1,\\
 z_1 &= m(m^8 + 2m^4 + 10),\quad  &z_2& = (m^4 + 3m^2 - 2)(m^4 - 3m^2 - 2),
\end{aligned}
\label{sol1m}
\end{equation}
where $m$ is an arbitrary parameter.

As  numerical examples, taking  $m=1$ and $m=2$, yields the solutions listed in Table I at S. No. 1 and S. No. 8  respectively. We could even assign rational values to $m$ in the solution \eqref{sol1m} and obtain integer solutions by appropriate scaling.

The parametric solution of Eq.~\eqref{quarteqn}, obtained from the  point $2P$ on the curve \eqref{ecweier},  may be written, in terms of the arbitrary parameter $m$, as follows:
\begin{equation}
\begin{aligned}
x_1 &= 12m^{21} - 282m^{17} + 1830m^{13} - 2256m^9 - 984m^5 + 480m,\\
 x_2 &= m^{24} - 240m^{16} + 1652m^{12} - 1800m^8 + 1668m^4 + 256, \\
y_1 &= m^{25} - 12m^{21} + 42m^{17} - 178m^{13} + 456m^9 + 2652m^5 - 224m,\\
 y_2 &= m^{24} - 6m^{20} - 99m^{16} + 737m^{12} - 672m^8 + 2160m^4 + 16,\\
 z_1 &= m^{49} - 12m^{45} - 126m^{41} + 970m^{37} + 30474m^{33}\\
 & \quad \quad - 405108m^{29}+ 1799754m^{25} - 3341016m^{21} + 3936600m^{17}\\
 & \quad \quad  - 1330304m^{13} + 4344720m^9 - 167040m^5 + 57856m,\\
 z_2 &= m^{48} - 78m^{44} + 1749m^{40} - 17069m^{36} + 79200m^{32}\\
  & \quad \quad - 183456m^{28} + 505284m^{24} - 3071484m^{20} + 10049220m^{16}\\
  & \quad \quad - 5711360m^{12} - 791232m^8 - 831552m^4 + 4096.
\end{aligned}
\end{equation}

The parametric solution corresponding to the point $3P$ gives the values of  $z_1$ and $z_2$ in terms of  polynomials of degree $\geq 120$ and is accordingly omitted.

We may also be able to find other points on the elliptic curve \eqref{ecweier} and thus obtain additional parametric solutions of Eq.~\eqref{quarteqn}. As an example, we found a point on the elliptic curve \eqref{ecweier} given by 
\[
\begin{aligned}
X &=4m^4(m^8 - 4m^4 - 14)^2/\{(m^4 + 3m^2 - 2)^2(m^4 - 3m^2 - 2)^2\}, \\
Y& = 36m^4(m^8 - 4m^4 - 14)(m^{16} - 11m^{12} + 30m^8 - 74m^4 - 8)\\
& \quad \quad \times \{(m^4 + 3m^2 - 2)(m^4 - 3m^2 - 2)\}^{-3},
\end{aligned}
\]
 and thus obtained the following  solution of Eq.~\eqref{quarteqn}:
\begin{equation}
\begin{aligned}
x_1 &= 6m^9 - 78m^5 + 24m, \\
x_2 &= m^{12} - 12m^8 + 72m^4 + 4,\\
 y_1 &= m^{13} - 6m^9 - 6m^5 + 28m,\\
 y_2 &= m^{12} - 9m^8 + 33m^4 + 16,\\
 z_1 &= m^{25} - 18m^{21} + 156m^{17} - 796m^{13} + 2394m^9 + 120m^5 + 400m, \\
z_2 &= m^{24} - 39m^{20} + 357m^{16} - 808m^{12} + 132m^8 - 2244m^4 + 64,
\end{aligned}
\end{equation}
where, as before, $m$ is an arbitrary parameter. 

\subsection{Solutions obtained by using  Pell's equation}
Several solutions listed in Table I cannot be generated by the parametric solutions obtained above. We also observed that several solutions in Table I had $x_1=1$. We therefore tried some experimentation and noticed that if we write,
\begin{equation}
\begin{aligned}
x_1 &= 1,\quad  &x_2 & = 2v,\\
y_1 &= 4v^2 + 1,\quad  &y_2 &= 2v(2v^2 + 1),\\
z_1 &= 4uv^2,\quad  & z_2&=8v^4 + 4v^2 + 1,
\end{aligned}
\label{subsPell}
\end{equation}
then Eq.~\eqref{quarteqn} reduces to
\begin{equation}
256v^8(u^2 + 3v^2 + 1)(u^2 - 3v^2 - 1)=0,
\end{equation}
and we can obtain integer solutions of \eqref{quarteqn} by choosing integers $u$ and $v$ satisfying the Pell's equation $u^2 - 3v^2=1$.  The two smallest nontrivial integer solutions of this equation namely $(u,\,v)=(2,\,1)$ and $(u,\,v)=(7,\,4)$ yield the two solutions of \eqref{quarteqn} listed in Table I at S. No. 1 and S. No. 5 respectively.

As already noted, rational solutions of Eq.~\eqref{quarteqn} immediately yield solutions in integers. Therefore  we need not restrict ourselves to choosing only integers $u$ and $v$ satisfying  the equation $u^2 - 3v^2=1$. A rational solution of this equation is readily seen to be as follows:
\begin{equation}
u = (t^2 + 3)/(t^2 - 3),\quad  v = 2t/(t^2 - 3),
\end{equation}
where $t$ is an arbitrary parameter. Now, using \eqref{subsPell}, we obtain the following solution of Eq.~\eqref{quarteqn}:
\begin{equation}
\begin{aligned}
x_1&=t^2-3,\quad & x_2&=4t, \\
y_1 &= (t^2 - 3)(t^2 + 9)(t^2 + 1), & y_2 &= 4t(t^2 + 2t + 3)(t^2 - 2t + 3),\\
z_1 &=16t^2(t^2 - 3)(t^2 + 3), & z_2 &= t^8 + 4t^6 + 86t^4 + 36t^2 + 81, \\
\end{aligned}
\end{equation}
where $t$ is an arbitrary parameter.

\section{Some open problems}
While we have seen in this paper that  parametric solutions of Eq.~\eqref{geneqn} exist when $n=3$ and  $n=4$, it would be of interest to find solutions of \eqref{geneqn} when $n > 4$. It is unlikely that integer solutions of \eqref{geneqn} exist with $n > 6$. Accordingly, it may be interesting to study the diophantine equations,
\begin{equation}
\sum_{i=1}^mx_i^n \sum_{i=1}^my_i^n =\sum_{i=1}^mz_i^n, \label{geneqnmterms}
\end{equation}
where $m$ is an integer as small as possible,  and $n \geq 5$.

More generally, it would be of interest to find  analogues of Brahmagupta's identity for cubes and higher powers, that is, we may seek identities of type \eqref{geneqnmterms}  where $n \geq 3$ while $x_i,\,y_i,\;i=1,\,2,\,\ldots,\,m$,  are arbitrary integers and $z_i,\;i=1,\,2,\,\ldots,\,m$ are integers determined by the integers $x_i,\,y_i$ as in the case of Brahmagupta's identity. 

We leave these problems for future investigations.

\noindent Ajai Choudhry, 13/4 A Clay Square, Lucknow - 226001, India

\noindent E-mail address: ajaic203@yahoo.com

\medskip 

\noindent Iliya Bluskov, P.O. Box 33031, West Vancouver, B.C. V7V 4W7, Canada

\noindent E-mail address: lotbook@telus.net

\medskip

\noindent Alexander James,  P.O.  Box 2004, Burwood North, NSW 2134, Australia

\noindent E-mail address: alexjamessolicitor@hotmail.com

\end{document}